\documentclass[10pt,a4paper]{article}

{}
{}
{}
\usepackage{color}
\usepackage{float}
\usepackage{hyperref}
\usepackage{multirow}
\usepackage{makecell}
\usepackage{float}
\usepackage{setspace}
\usepackage{array}
\usepackage{supertabular}
\usepackage{multicol}
\usepackage{geometry}
\usepackage{cite} 
\usepackage{mathtools}
\usepackage{amssymb}
\usepackage{graphicx}
\usepackage{epsfig}
\usepackage[ruled,linesnumbered]{algorithm2e}
\usepackage{comment}

\begin{document}
 \begin{center}
{\LARGE\textbf{Improving Operational Feasibility of Low-voltage}} \\
\vspace{3pt}
{\LARGE\textbf{Distribution Network by Phase-switching Devices}} \\
\vspace{5pt}
{ {Bin Liu$^{1*}$, Ke Meng$^{1}$, Peter K.C. Wong$^{2}$, Zhao Yang Dong$^{1}$, Cuo Zhang$^{1}$, Bo Wang$^{1}$, Tian Ting$^{3}$, Qu Qi$^{4}$}}\\
\vspace{10pt}
{\small $^{1}$School of Electrical Engineering and Telecommunications, The University of New South Wales, Sydney 2052, Australia\\
$^{2}$Jemena Electricity Networks (VIC) Ltd., 567 Collins Street, Melbourne 3000, Australia\\
$^{3}$Ausnet Services, 2 Southbank Boulevard, Melbourne 3006, Australia\\
$^{4}$State Grid International Development Pty. Ltd., 88 Chang’an Avenue, Beijing 100031, China\\
$^{*}$bin.liu@unsw.edu.au, eeliubin@hotmail.com}\\
\vspace{10pt}
\textbf{Keywors}: Linearisation technique, mixed-integer constraints, phase-switching device, residential PV, unbalanced power flow 
\end{center}

\setlength{\columnsep}{10pt}  
\begin{multicols}{2}

\section*{Abstract}
High penetration of residential photovoltaic (PV) in low-voltage distribution networks (LVDN) makes the unbalance issue more sever due to the asymmetry of generation/load characteristic in different phases, and may lead to infeasible operation of the whole network. Phase switching device (PSD), which can switch the connected phase of a residential load as required, is viable and efficient equipment to help address this issue. This paper, based on three-phase power flow (TUPF) formulation, aims to investigate the benefit of PSD on improving the operation feasibility of LVDN. The linear model of TUPF with PSD is presented to effectively take the flexible device into account, which can be conveniently used to seek the PSD positions that lead to a viable operation strategy of the original problem when infeasibility is reported by the traditional iteration-based algorithm. Case study based on a practical LVDN in Australia demonstrates the efficacy of the proposed method.

\section{Introduction}
To reduce the consumption of irreversible fossil fuels and alleviate the concerns of climate change issue, renewable energy development has been promoted by most countries world-widely. Among various renewable power generation methods, wind and photovoltaic (PV) generation are most popular in practice due to their technology maturity and sustained decrease in cost. With encouraging polices and incentives from government, Australia is experiencing a remarkable renewable energy development, taking the top spot worldwide in the penetration level of residential PV installation in low-voltage distribution network (LVDN) . The annual report released by Australian Energy Council shows that it was another record-breaking year for PV development in Australia in 2018, with the installed residential capacity reaching over 1.4 GW, which increases by 20\% compared with 2017. By the end of 2018, cumulative installed capacity of residential PV in Australia stood at 7.98 GW with more than 2 million installations across the nation, and the numbers continues to grow  \cite{RN192}. 

Power utilities in Australia run extensive four-wire (230/400 V) grid along the streetscape as in many other countries. Multiple points in the LVDN and the neutral conductor at the distribution transformer (DT) all are earthed \cite{RN61}. Most residential customers are powered by a bundled cable, which consists of a powered phase and a neutral phase from the nearest pole. Although increasing residential PV generations significantly reduce the bills of electricity customers, it also makes the existing unbalance issue more severe in LVDN. The unbalance issue, which is caused by unbalanced load profile traditionally, can be worsened because PV panels are installed out of customers' wills without careful design beforehand \cite{RN52}. Unbalances in LVDN can cause lots of operational problems, e.g. high level of neutral current that may lead to protection system malfunction and increased power loss, low electricity supply quality that can shorten the lives of both electrical equipment in the network and appliances of residential customers, and significant voltage drop or heavy loaded feeders in one phase that decrease the usage efficiency of the whole network \cite{RN103,RN104}. From operator's perspective, the unbalance issue can undermine the operational feasibility of LVDN, i.e. leading to infeasible operation by some generation/load profiles. Reflected in power system analysis, infeasible operational situation leads to insolvable case of three-phase unbalanced power flow (TUPF).  

To mitigate the unbalance issue, one cost-efficient choice is to employ phase-switching device (PSD) to change the phase that each residential customer is connected depending on the operational state of the LVDN \cite{RN61,RN44}. As deploying PSD provides more options to operate the LVDN, the corresponding TUPF may become solvable due to the enlarged feasible region, thus improving the operational feasibility. 

On the topic of solving TUPF, most researchers mainly focus on iteration-based method, e.g. the rigid approach considering the sparsity of distribution network \cite{RN25}, the fast decoupled method \cite{RN26}, the newton-raphson method \cite{RN31}, the current-injection method \cite{RN69}, the direct approach based on bus-injection to branch-current matrix and branch-current to bus-voltage matrix \cite{RN37}, and the approach based on complex theory in $\alpha\beta 0$ stationary reference frame \cite{RN70}. As iteration-based algorithms relies on the initialized points and the nonlinear formulation of TUPF is not convenient when applied in optimization problem, recent researches focus on investigating the convergence of iteration-based algorithms \cite{RN55,RN59,RN23,RN22} and the linearization of TUPF \cite{RN38,RN39,RN29,RN59}. However, few literature has been reported in solving TUPF of LVDN equipped with PSDs. This is because equipping PSDs in LVDN will introduce integer variables to the TUPF model, thus making the iteration-based algorithm much more complicated. 

Based on reported methods in handling non-convex parts in TUPF model, this paper presents a mixed-integer linear (MIL) formulation of TUPF for LVDN equipped with PSDs, which enables us to investigate the operational feasibility improvement brought by PSD. The efficacy is supported by accurate iteration-based algorithm to solve TUPF and demonstrated on a practical system in Australia. 
The remaining of the paper is organized as follows. The mathematical formulation of TUPF for LVDN with PSDs is presented in Section II, followed by the non-convexity analysis and solution technique in Section III. Case study based on a practical system in Australia will be performed in Part IV. The paper is concluded in Part V.
\section{TUPF for LVDN with PSDs}
\begin{figure}[H]
	\centering\includegraphics[scale=0.37]{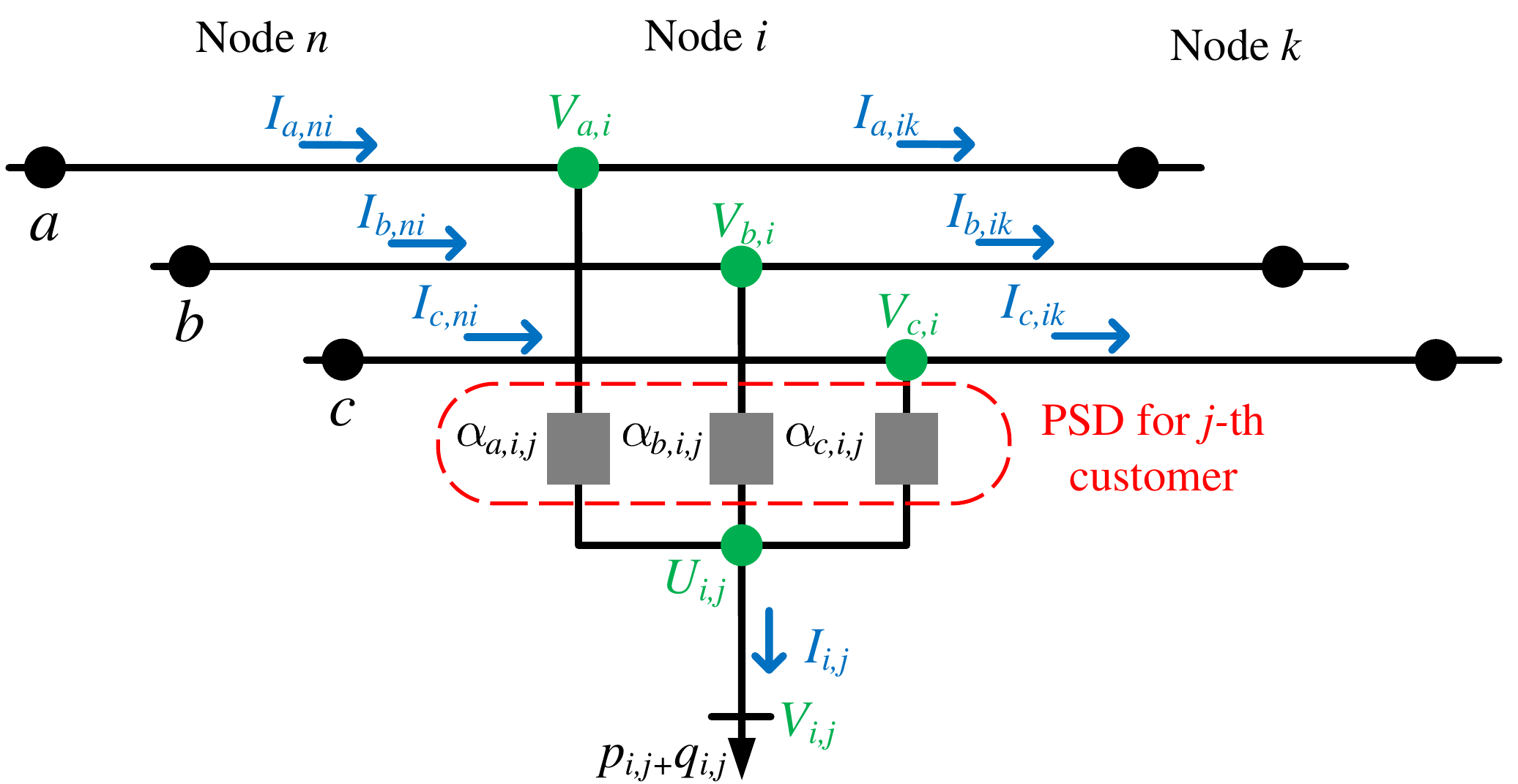}
	\caption{Illustration of operating PSD for flexible customer}
	\label{fig-psd}
\end{figure}

The formulation are based on Fig.\ref{fig-psd}.
In the formulated problem, $\mathcal {C}_i$ is used to represent the set of customers connected to node $i$, $\mathcal{F}_i$ and $\mathcal{X}_i$ as flexible customers (with PSD installed) and inflexible customers (without PSD), respectively. Obviously, we have $\mathcal{F}_i\cap \mathcal{X}_i=\emptyset$ and $\mathcal{F}_i\cup \mathcal{X}_i=\mathcal{C}_i$. Moreover, $V,U$ will represent voltage with $X$ and $Y$ being the real and image part of $V$, respectively; $I$ represent current with $J$ and $W$ being the real and image parts, respectively; The subscript $\phi/\psi,i,j,ik$ represent phase $\phi/\psi$ that belongs to $\{a,b,c\}$, node $i$, the $j^\text{th}$ customer at node $i$ and line $ik$, respectively. Other parameters or variables will be explained right after their appearances.

Moreover, we have the following assumptions for formulating TUPF. 
\begin{enumerate}
	\item The root node, i.e. the low-voltage side of DT, is selected as the balancing node. Therefore, its voltage is a known parameter.
	\item All residential customers are $PQ$ bus and powered by single-phase. In other words, the active/reactive powers of all customers are known parameters.
	\item Lines between any two poles in LVDN are constructed with four-wire (phase $a,b,c$ and zero earthed conductor), which is the general case in Australia \cite{RN61}.
\end{enumerate}

In the TUPF formulation, the Ohm's law for each line in main feeder must be satisfied as follows.
\begin{eqnarray}
\label{ol-1}
V_{\phi,i}-V_{\phi,k}=\sum_\psi{Z_{i,k}^{\phi\psi}I_{\phi,ik}}~\forall \phi,\forall ik
\end{eqnarray}
where 
$Z_{i,k}^{\phi\psi}$ is the mutual impedance between phase $\phi$ and $\psi$ of line $ik$.

The KCL must be satisfied to ensure the current balance at each node as shown in \eqref{kcl-1}. Moreover, the current injected to the node from its connected customers should be consistent as shown in \eqref{kcl-2} and \eqref{kcl-3}, and each customer can only be connected to one phase as indicated by \eqref{kcl-4} and \eqref{kcl-5}.
\begin{subequations}
	\begin{eqnarray}
	\label{kcl-1}
	\sum_{n:n\rightarrow i}{I_{\phi,ni}}-\sum_{k:i\rightarrow k}{I_{\phi,ik}}=\sum_{j\in \mathcal {C}_i}{I_{\phi,i,j}}~\forall \phi,\forall i\neq x\\
	\label{kcl-2}
	I_{\phi,i,j}=\alpha_{\phi,i,j}I_{i,j}~\forall \phi,\forall i,\forall j\in \mathcal{F}_i\\
	\label{kcl-3}
	I_{\phi,i,j}=\mu_{\phi,i,j}I_{i,j}~\forall \phi,\forall i,\forall j\in \mathcal{X}_i\\
	\label{kcl-4}
	\alpha_{\phi,i,j}\in \{0,1\}~\forall \phi,\forall i,\forall j\\
	\label{kcl-5}
	\sum_\phi{\alpha_{\phi,i,j}=1}~\forall i,\forall j
	\end{eqnarray}
\end{subequations}
where
$x$ is the index of root node; 
$\alpha_{\phi,i,j}$ is binary variable indicating whether the $j^\text{th}$ customer is connected to phase $\phi$ of node $i$; 
$\mu_{\phi,i,j}$ is known parameter indicating the initial phase-position of $j^\text{th}$ customer.

Ohm's law should be satisfied for each connecting line and the terminal voltages of PSD should be identical, leading to
\begin{subequations}
	\begin{eqnarray}
	\label{kcl-6}
	U_{i,j}-V_{i,j}=Z_{i,j}I_{i,j}~\forall i,\forall j\in\mathcal{C}_i\\
	\label{kcl-7}
	U_{i,j}=\sum_{\phi}{\alpha_{\phi,i,j}V_{\phi,i}}~\forall i,\forall j\in\mathcal{F}_i\\
	\label{kcl-8}
	U_{i,j}=\sum_{\phi}{\mu_{\phi,i,j}V_{\phi,i}}~\forall i,\forall j\in\mathcal{X}_i
	\end{eqnarray}
\end{subequations}
where 
$Z_{i,j}$ is the impedance of the line from node $i$ to its $j^\text{th}$ customer. 

For each customer, the power balance equations are expressed as	
\begin{eqnarray}
\label{pbeq}
\label{pb-1}
I_{i,j}=\frac{P_{i,j}-jQ_{i,j}}{V_{i,j}^*}~\forall i,\forall j\in\mathcal{C}_i
\end{eqnarray}
where $P_{i,j},Q_{i,j}$ are net active and reactive demands of $j^\text{th}$ customer at node $i$, respectively. 

Particularly, the voltage of root node is assumed to be known, i.e.
\begin{eqnarray}
\label{vroot}
V_{\phi,x}=V_{\phi}^0~\forall \phi
\end{eqnarray}
where $V_{\phi,s}^0$ is the known voltage of root node at phase $\phi$.

When the phase positions of residential customers are known, the TUPF formulation is equivalent to that expressed by bus-injection-to-bus-current (BIBC) and branch-current-to-bus-voltage (BCBV) matrices in \cite{RN37}, thus can be solved by the iteration-based algorithm proposed there. However, when with PSDs, the algorithm is no longer viable due to introduced integer variables, which motivate us to develop other efficient method to solve the challenging problem in the next section.

\section{Solution Technique}
Noting that the non-convex parts in TUPF formulation is due to the bilinear terms in \eqref{kcl-2} and \eqref{kcl-7}, and the division operator in \eqref{pb-1}. The main idea in this paper is to reformulate or lineairze the non-convex parts to make it efficiently solvable, as we explain next. 
\subsection{Linearizing bilinear terms}
Noting that the bilinear terms in \eqref{kcl-2} and \eqref{kcl-7} are the product of a binary and continuous variables, they can be exactly reformulated as mixed-integer linear (MIL) constraints noting the following equivalent formulation \cite{RN30,BLiu-ITES}.
\begin{eqnarray}
\label{RLP-1}
\left.\begin{array}{r}
z=xy\\
x\in \{0,1\}\\
y\in[y^\text{min},y^\text{max}]
\end{array}
\right\}
\Leftrightarrow 
\left\{\begin{array}{r}
xy^\text{min}\le z \le xy^\text{max}\\
(x-1)y^\text{max}\le z-y\\
\le(x-1)y^\text{min}
\end{array}
\right.\
\end{eqnarray}

Based on \eqref{RLP-1}, \eqref{kcl-2} can be exactly reformulated as
\begin{eqnarray}
\label{kcl-2-mil}
I_{\phi,i,j}=z^J_{\phi,i,j}+jz^W_{\phi,i,j}\\
\alpha_{\phi,i,j}J^\text{min}\le z^J_{\phi,i,j}\le\alpha_{\phi,i,j}J^\text{max}\\
(\alpha_{\phi,i,j}-1)J^\text{max}\le z^J_{\phi,i,j}-J_{i,j}\le(\alpha_{\phi,i,j}-1)J^\text{min}\\
\alpha_{\phi,i,j}W^\text{min}\le z^W_{\phi,i,j}\le\alpha_{\phi,i,j}W^\text{max}\\
(\alpha_{\phi,i,j}-1)W^\text{max}\le z^W_{\phi,i,j}-W_{i,j}\le(\alpha_{\phi,i,j}-1)W^\text{min}
\end{eqnarray}

Similarly, \eqref{kcl-7} can be reformulated as
\begin{eqnarray}
\label{kcl-7-mil}
U_{i,j}=\sum_{\phi}{z^E_{\phi,i,j}}+j\sum_{\phi}{z^F_{\phi,i,j}}\\
\alpha_{\phi,i,j}X^\text{min}\le	z^E_{\phi,i,j}\le \alpha_{\phi,i,j}X^\text{max}\\
(\alpha_{\phi,i,j}-1)X^\text{max}\le z^E_{\phi,i,j}-X_{\phi,i}\le(\alpha_{\phi,i,j}-1)X^\text{min}\\
\alpha_{\phi,i,j}Y^\text{min}\le	z^F_{\phi,i,j}\le \alpha_{\phi,i,j}Y^\text{max}\\
(\alpha_{\phi,i,j}-1)Y^\text{max}\le z^F_{\phi,i,j}-Y_{\phi,i}\le(\alpha_{\phi,i,j}-1)Y^\text{min}
\end{eqnarray}

Thus, the con-convex bilinear terms in \eqref{kcl-2} and \eqref{kcl-7} are exactly reformulated as MIL constraints. 

\subsection{Linearizing power balance equation}
The linearization of power balance equation is mainly based on the assumption that the voltage angles of all nodes in each phase are with limited difference, which has been demonstrated in \cite{RN38,RN50,RN66,RN30,RN67}. To make this more clear in formulation, we assume VM limit at node $i$ is $[V^\text{min}_{i},V^\text{max}_{i}]$ for any phase, and VA limit in phase $\phi$ for any node is $[\delta^\text{min}_\phi,\delta^\text{max}_\phi]$. Further, VA is assumed to be centered at $\delta_\phi$, which leads to $\delta^\text{min}_\phi=\delta_\phi-\Delta \delta,\delta^\text{max}_\phi=\delta_\phi+\Delta\delta$.

Specifically, the linearziation idea is firstly approximating $1/V^*$ in \eqref{pb-1} by a linear function, and then to further linearize the power balance equations. 
\subsubsection{Linearizing $1/V^*$}
For $V_{i,j}$ in \eqref{pb-1}, we here take the following generalized form to illustrate how to approximately linearize it.
\begin{eqnarray}
\label{lp-v}
f(V)=\frac{1}{V^*}=\frac{1}{X-jY}=\frac{X+jY}{|V|^2}
\end{eqnarray} 
where $V=X+jY$.

To linearize \eqref{lp-v}, one method is based on linear expansion of a function in complex domain as discussed in \cite{RN29}. For $f(V)$, it can be approximated by $2-V^*$ when $V^*$ is around $1+j0$. For $V^*$ with VA centered at $\delta$, the following linearization can be derived.
\begin{eqnarray}
f(V)=\frac{1}{V^*}\approx 2e^{j\delta}-V^*e^{j2\delta}
\end{eqnarray} 

We denote this method as complex-based method (CBM) throughout the context.

Another method is employing least square method (LSM) to linearize $f_X(V)=X/|V|^2$ and $f_Y(V)=Y/|V|^2$, by the following expressions, respectively \cite{RN39,BLiu-Bcoff}.
\begin{subequations}
	\label{fxy}
	\begin{eqnarray}
	f_X(V)\approx k^XX+k^YY+b^X\\
	f_Y(V)\approx h^XX+h^YY+b^Y
	\end{eqnarray} 
\end{subequations} 
where $k^X,k^Y,b^X,h^X,h^Y,b^Y$ are parameters to be fitted.

To fit all parameters in \eqref{fxy}, a sufficient number of known points must be provided, which can be sampled in the specified region as follows taking the phase $\phi$ of node $i$ as an example.
\begin{itemize}
	\item Discretise the feasible region of voltage magnitude (VM) of node $i$, i.e. $[V_i^\text{min},V_i^\text{max}]$, as $M$ points, the set of which is denoted as $\mathcal{VM}=\{V_i^1,\cdots,V_i^M\}$. 
	\item Discretise the feasible region of voltage angle (VA) of phase $\phi$, i.e. $\delta_\phi^\text{min},\delta_\phi^\text{max}]$, as $N$ points, the set of which is denoted as $\mathcal{VA}=\{\delta_\phi^1,\cdots,\delta_\phi^N\}$.
	\item For each combination of elements in $\mathcal{VM}$ and $\mathcal{VA}$, say $(V_i^m,\delta_\phi^n)$, exact values of $X_\phi^{m,n},Y_\phi^{m,n},f^{m,n}_{X,\phi},f^{m,n}_{Y,\phi}$ are calculated as
	\begin{eqnarray}
		X_\phi^{m,n}=V_i^m\cos{\delta_\phi^n},~Y_\phi^{m,n}=V_i^m\sin{\delta_\phi^n}\nonumber\\
		f_{X,\phi}^{m,n}=\frac{X_\phi^{m,n}}{(X_\phi^{m,n})^2+(Y_\phi^{m,n})^2}\nonumber\\
		f_{Y,\phi}^{m,n}=\frac{Y_\phi^{m,n}}{(X_\phi^{m,n})^2+(Y_\phi^{m,n})^2}\nonumber
	\end{eqnarray} 
	\item All the sampled points, i.e.
	$$(X_\phi^{m,n},Y_\phi^{m,n},f^{m,n}_{X,\phi},f^{m,n}_{Y,\phi})~(\forall m,\forall n)$$
	will be used to determine the parameters in \eqref{fxy}. 
\end{itemize}

\begin{figure}[H]
	\centering\includegraphics[scale=0.19]{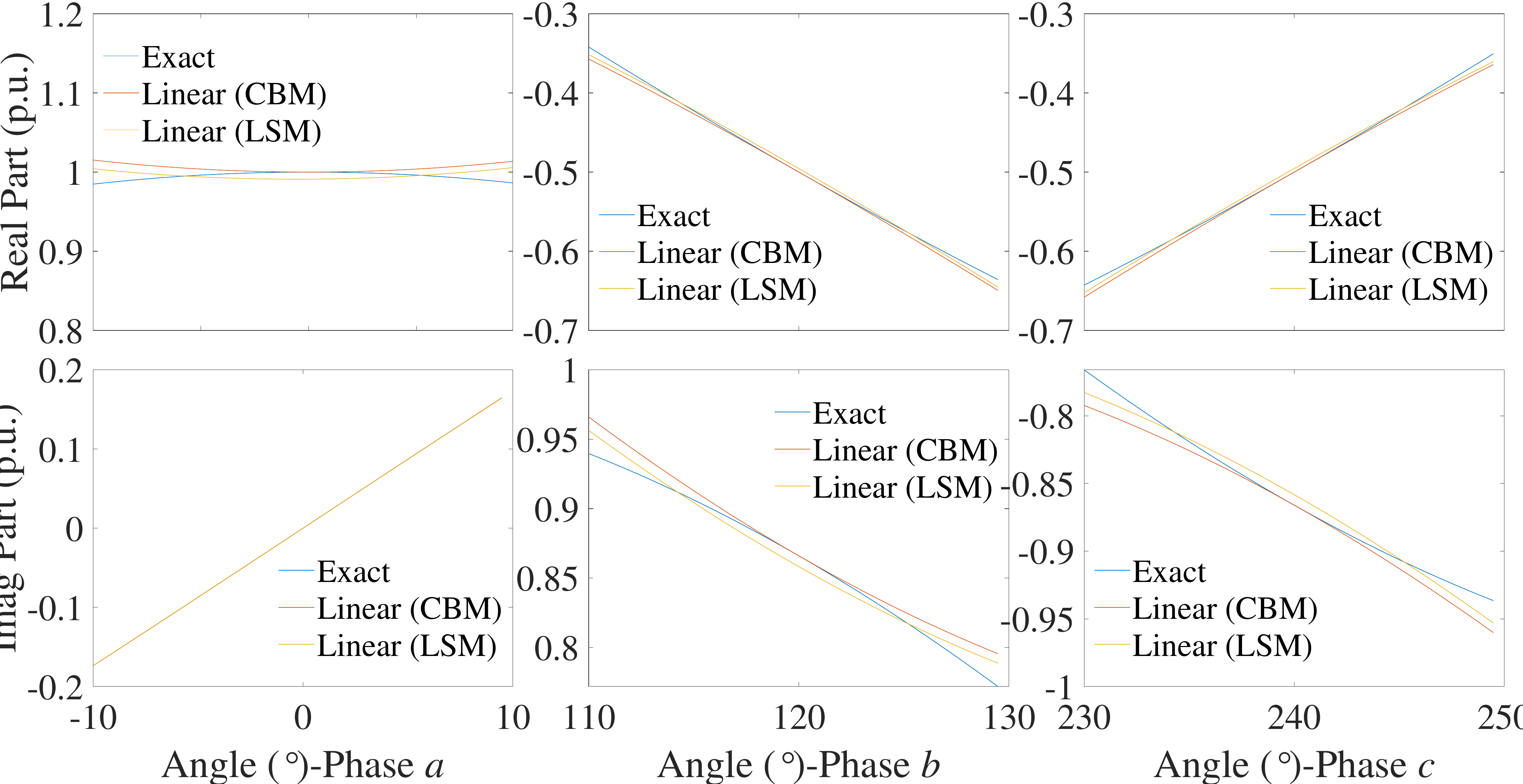}
	\caption{Comparison of different methods to linearize $f_X(V)$ and $f_Y(V)$ (Only the fitting results when $|V|=1$ is presented with respect to different VAs)}
	\label{fig-linearapprox}
\end{figure}

Moreover, parameters in \eqref{fxy} should be calculated for each phase separately. The analytical expression of the parameters can be found in \cite{RN38} and is omitted here for simplicity.

The comparison of the two methods are showed in Fig.\ref{fig-linearapprox}, where the both methods are based on $V^\text{min}=0.95,V^\text{max}=1.05$, and $\delta_a=0^\circ,\delta_b=120^\circ,\delta_c=240^\circ$, and $\Delta\delta=10^\circ$. The simulation results show both of the methods are with high accuracy. In this paper, the LSM will be employed for TUPF linearization.

\subsubsection{Linearizing power balance equation}
With \eqref{fxy}, \eqref{pb-1} can be reformulated as
\begin{eqnarray}
\label{pb-1-1}
\begin{split}
\label{pb-1-2}
I_{i,j}=(P_{i,j}-jQ_{i,j})[k^X_{i,j}X_{i,j}+k^Y_{i,j}Y_{i,j}+b^X_{i,j}\\
+j(h^X_{i,j}X_{i,j}+h^Y_{i,j}Y_{i,j}+b^Y_{i,j})]~\forall i,\forall j\\
k^X_{i,j}=\sum_\phi{\varkappa_{\phi,i,j}k^X_{\phi}},~h^X_{i,j}=\sum_\phi{\varkappa_{\phi,i,j}h^X_{\phi}}~\forall i,\forall j\\
k^Y_{i,j}=\sum_\phi{\varkappa_{\phi,i,j}k^Y_{\phi}},~h^Y_{i,j}=\sum_\phi{\varkappa_{\phi,i,j}h^Y_{\phi}}~\forall i,\forall j\\
b^X_{i,j}=\sum_\phi{\varkappa_{\phi,i,j}b^X_{\phi}},~b^Y_{i,j}=\sum_\phi{\varkappa_{\phi,i,j}b^Y_{\phi}}~\forall i,\forall j
\end{split}
\end{eqnarray}
where $\varkappa_{\phi,i,j}$ equals to $\mu_{\phi,i,j}$ for $j\in\mathcal{X}_i$ and $\alpha_{\phi,i,j}$ for $j\in\mathcal{F}_i$

Obviously, for inflexible customers, \eqref{pb-1-1} is linear. However, for flexible customers, expanding \eqref{pb-1-2} will lead to non-convex bilinear terms. Noting that each bilinear term is again the product of a integer variable and a continuous variable, it can be exactly linearized according to \eqref{RLP-1} as well.

\subsection{Algorithm}
Based on the solution techniques discussed previously, the algorithm to improve operational feasibility of the LVDN with PSD is given as Algorithm \ref{alg-1}.
\begin{algorithm}[H]
	\footnotesize  
	\label{alg-1}
	\caption{Algorithm to improve operational feasibility of TUPF in LVDN with PSD}
	\KwIn{Customer load and network parameters, initial positions of PSDs $\mu_{\phi,i,j}(\forall \phi, \forall i,\forall j\in\mathcal{C}_i)$}
	\KwResult{$V_{\phi,i}(\forall \phi,\forall i)$ and $\alpha_{\phi,i,j}(\forall \phi,\forall i,\forall j\in\mathcal{F}_i)$}
	Solve NTUPF with $\alpha_{\phi,i,j}=\mu_{\phi,i,j}(\forall \phi,\forall i,\forall j\in \mathcal{F}_i)$.\\
	\If {NTUPF is feasible}{End and report the solution.}
	\Else
	{Solve LTUPF regarding $\alpha_{\phi,i,j}$ as variables.\\
		\If{LTUPF is infeasible}{End and report infeasibility.}
		\Else{Solve NTUPF with $\mu_{\phi,i,j}$ fixed as $\alpha_{\phi,i,j}$ calculated from LTUPF.\\
			\If{NTUPF is feasible}{End and report the solution.}
			\Else{End and report infeasibility.}}}
\end{algorithm}

Several remarks on the formulation and solution technique are given below.
\begin{enumerate}
	\item The linear TUPF with PSD (LTUPF) is in MIP forms, where the error only arises from linearizing $1/V$ compared with its exact formulation, i.e. non-convex TUPF (NTUPF). The resulted errors in LTUPF by linearizing $1/V$ will be further tested in case study. 
	\item Due to introduced integer variables, even LTUPF may have multiple solutions. To demonstrate the benefit of deploying PSD more clear, we take the solution achieving best balance among three phases, i.e. minimum $U_x$ in \eqref{objx}, as the one to be investigated. 
	\begin{eqnarray}
		\label{objx}
		F=\min\left\{U_x\Bigg|
		\begin{array}{c}
		|P_\phi-P_\psi|\le U_x~\forall \phi,\psi\\
		|Q_\phi-Q_\psi|\le U_x~\forall \phi,\psi
		\end{array}\right\}
	\end{eqnarray}
	where $P_\phi=\sum_{i}(\sum_{j\in\mathcal{X}_i}{\mu_{\phi,i,j}P_{i,j}}+\sum_{j\in\mathcal{F}_i}{\alpha_{\phi,i,j}P_{i,j}})$ and 
	$Q_\phi=\sum_{i}(\sum_{j\in\mathcal{X}_i}{\mu_{\phi,i,j}Q_{i,j}}+\sum_{j\in\mathcal{F}_i}{\alpha_{\phi,i,j}Q_{i,j}})$.
	
	It is noteworthy that seeking the preferred LTUPF is different from solving the more complicated three-phase unbalanced optimal power flow (TUOPF) problem, where the VM limits and ampacities of conductors must be considered as well. LTUPF is solved by commercial solvers Gurobi 8.0.1 \cite{gurobi} in this paper.
	
	\item Similar to single-phase direct-current power flow (DCPF) and alternative-current power flow (ACPF), solving LTUPF may lead to actually infeasible solutions\footnote{The feasible solution here refers to a combination of net active/reactive power demands in LVDN.}. The feasibility is verified or judged by again running the NTUPF with updated power demand. However, infeasibility reported by NTUPF does not conclude that the solution reported by LTUPF is infeasible, because checking the feasibility of a combination of load profiles is NP-hard \cite{BLiu-Bcoff,BPF-2}. This is beyond the scope of this paper and will not be discussed in detail.
\end{enumerate}

\section{Case Study}
\subsection{Simulation setup}
\begin{figure}[H]
\centering\includegraphics[scale=0.41]{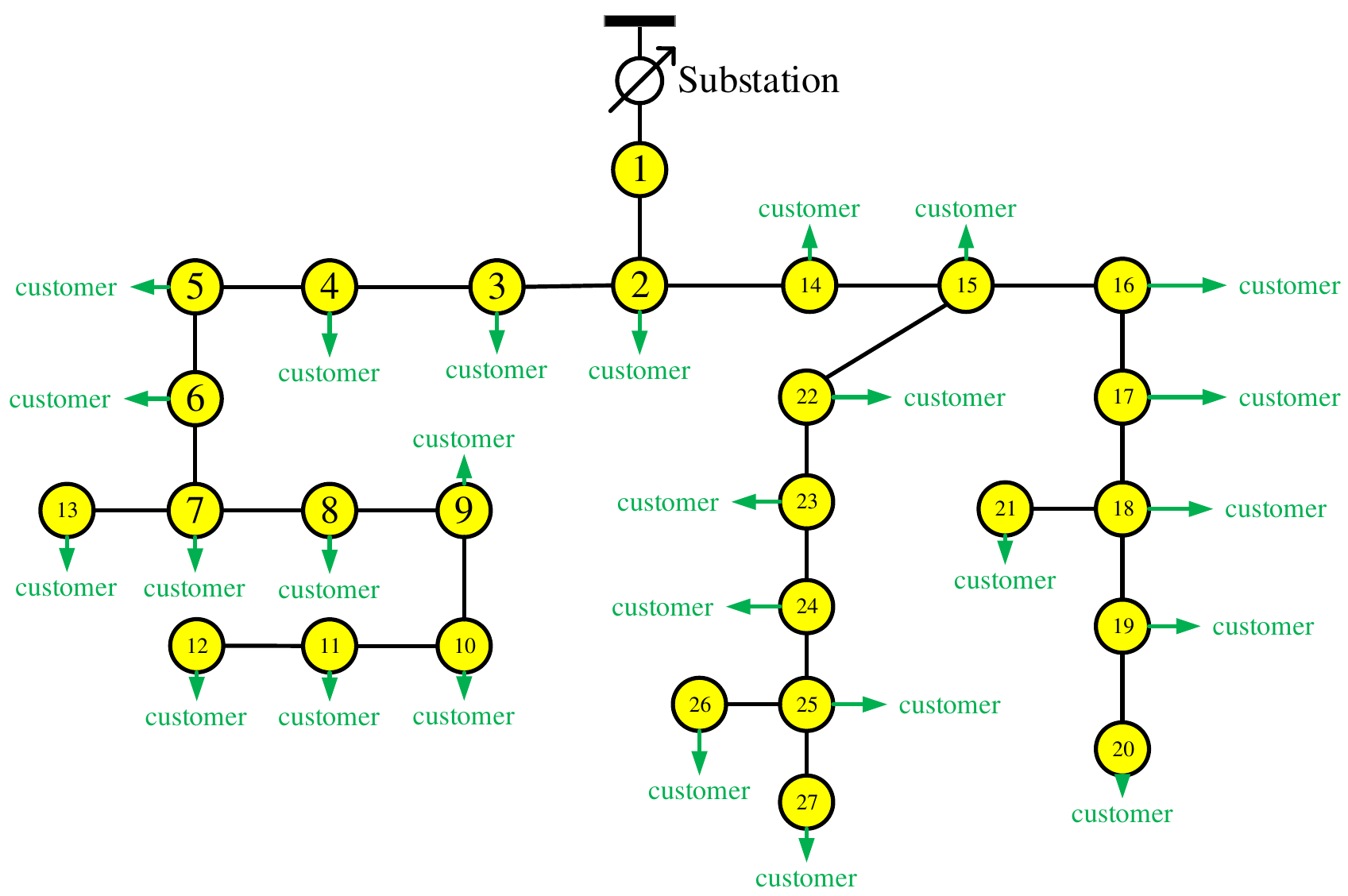}
\caption{Topology of the studied LVDN in Australia}
\label{fig-ausnet-paper}
\end{figure}

\begin{figure}[H]
	\centering\includegraphics[scale=0.17]{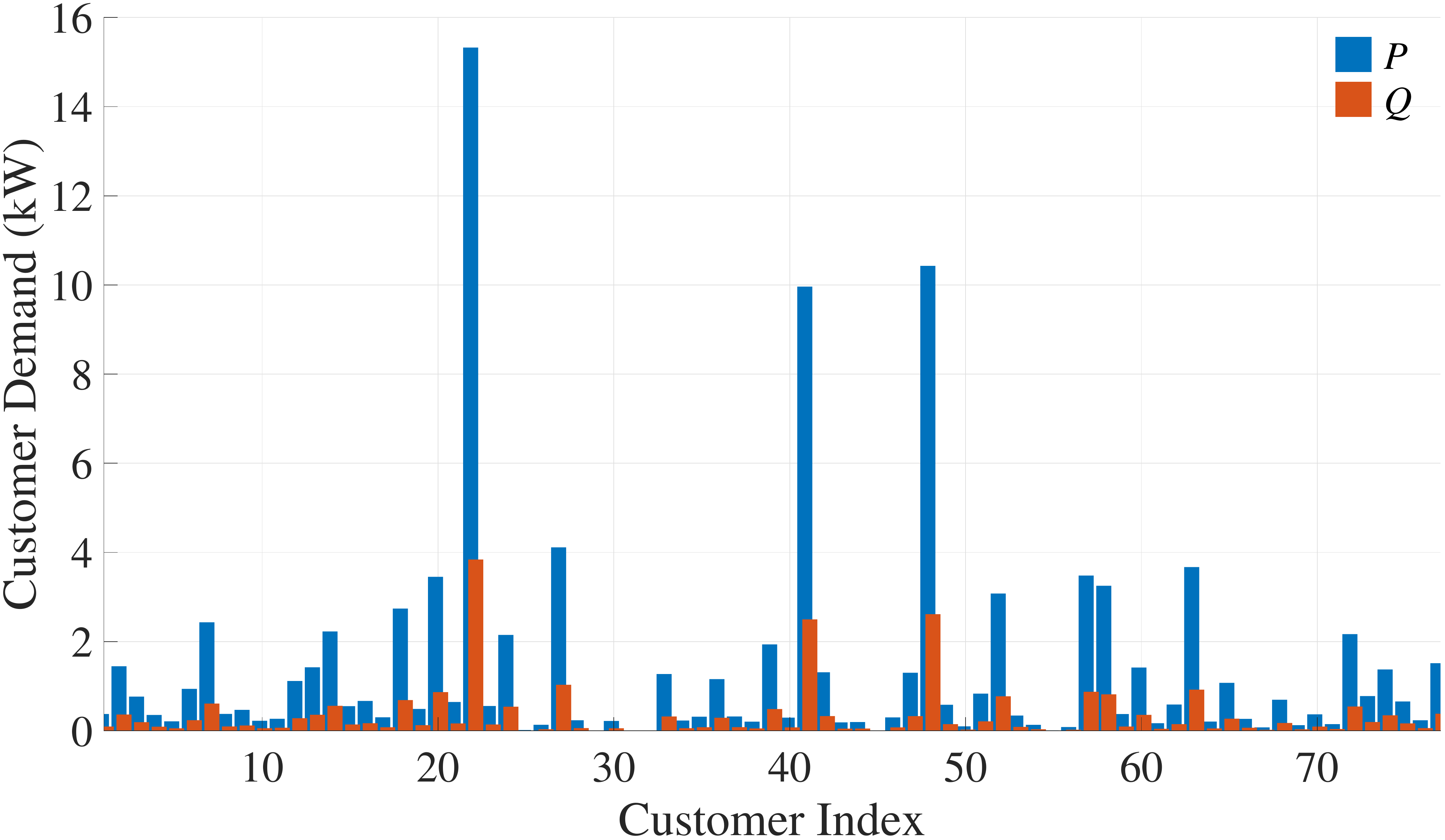}
	\caption{Active/reactive demand of residential customers}
	\label{fig-cusPQ}
\end{figure}

A practical LVDN in Australia is employed in this paper to investigate how the PSD help improve the operational feasibility. The topology of the system is as shown in Fig.\ref{fig-ausnet-paper}. In the LVDN, 77 single-phase powered residential customers are connected, where customer $1,9,17\cdots,65$ and $73$ are equipped with PSDs. The active/reactive demand of residential customers are shown in Fig.\ref{fig-cusPQ}.

In the formulation, the voltage of root node is known parameter, which is set as 
$$V^0=[1.05e^{j0},1.05e^{j\frac{2\pi}{3}},1.05e^{-j\frac{2\pi}{3}}]^T$$

Two cases based on the practical system will be studied, i.e. 
\begin{itemize}
	\item Case I: This case is totally based on the parameters given above and the purpose is to investigate the accuracy of LTUPF compared with NTUPF.
	\item Case II: Parameters in this case is the same as those in Case I except that the all customer load located at phase $c$ are doubled.
\end{itemize}

\subsection{Case I}
\begin{figure}[H]
	\centering\includegraphics[scale=0.17]{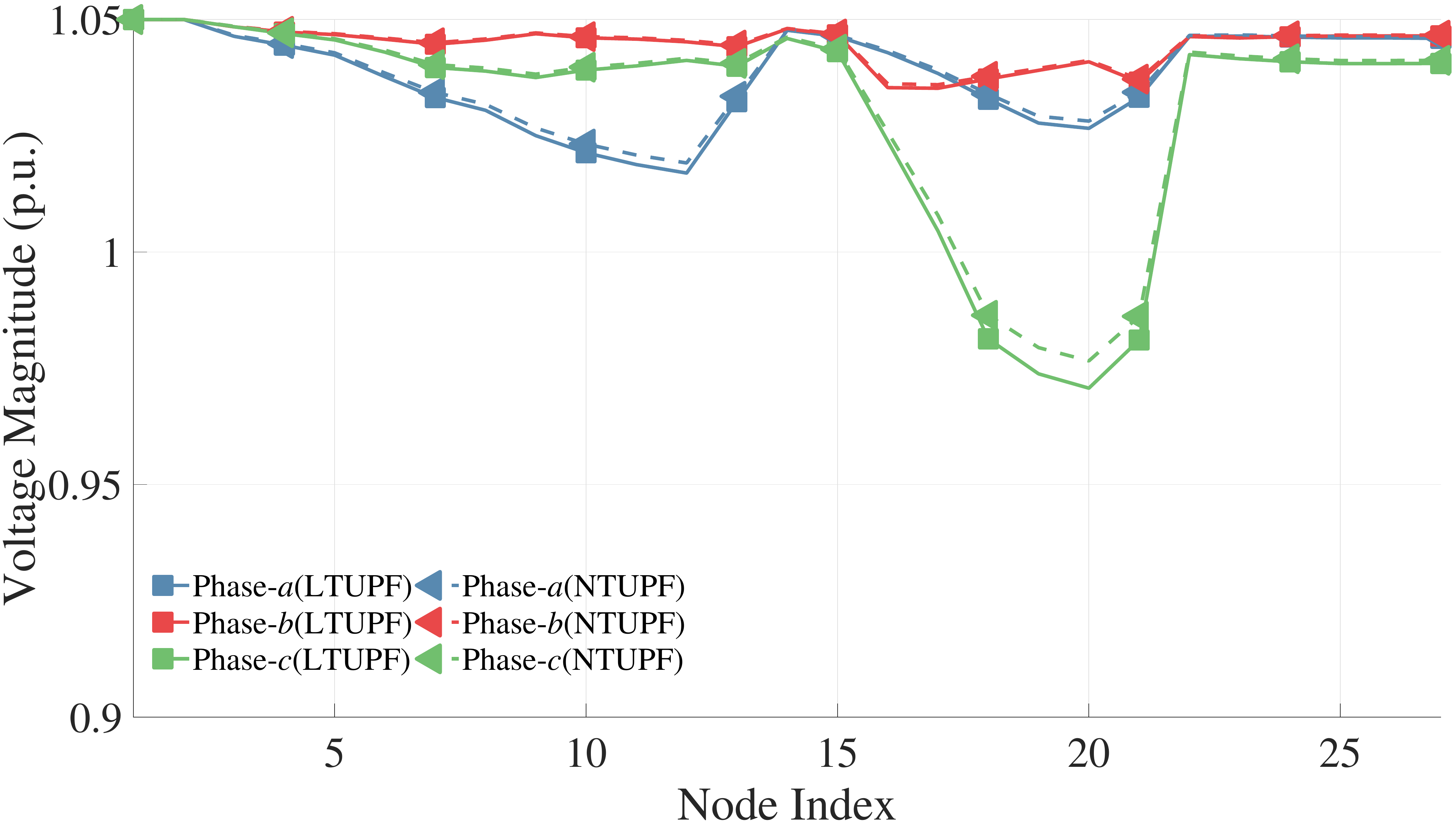}
	\caption{Calculated VMs for Case I}
	\label{fig-VMcvged}
\end{figure}

\begin{figure}[H]
	\centering\includegraphics[scale=0.17]{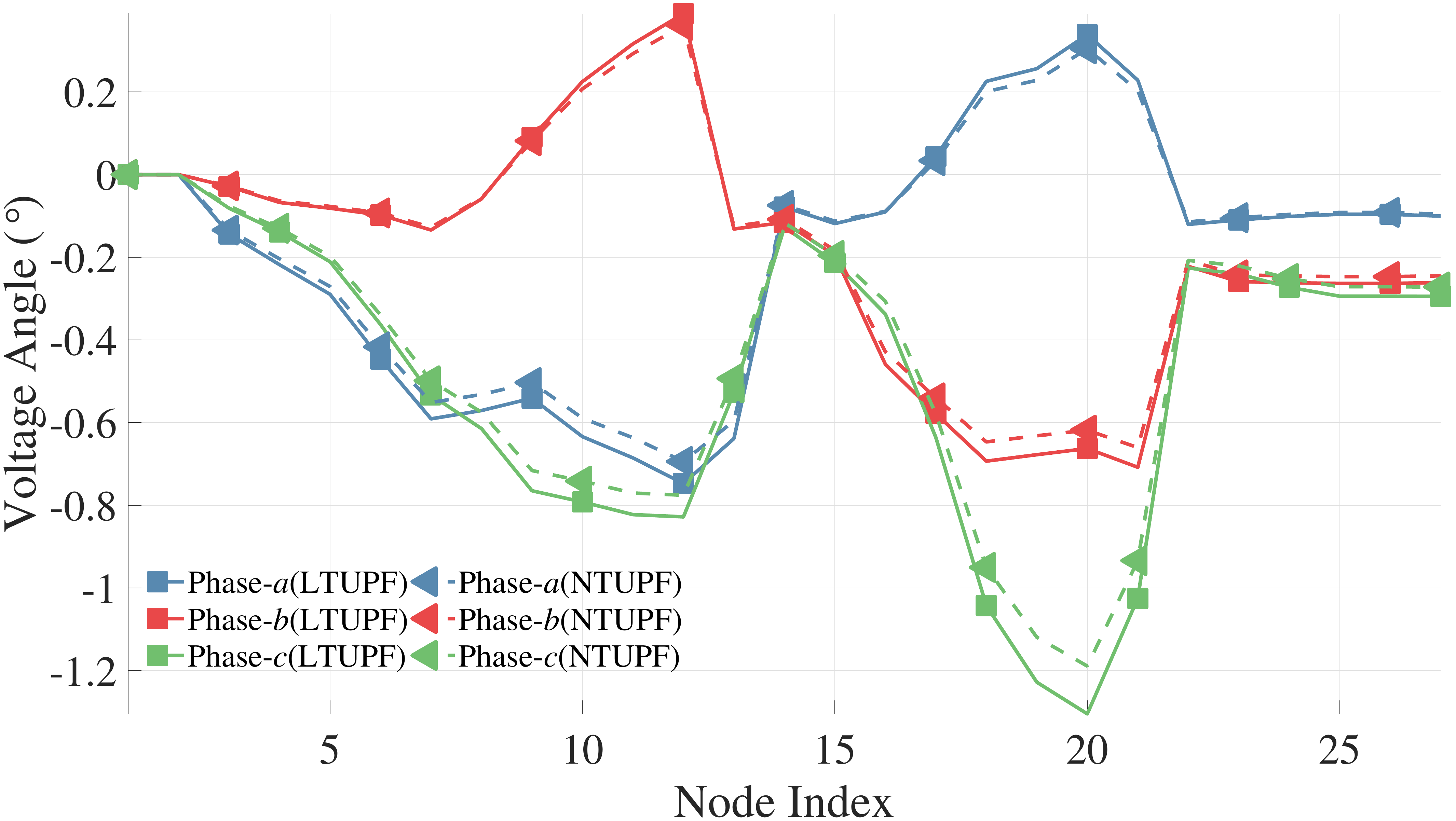}
	\caption{Calculated VAs for Case I (VAs for phase $b$ and $c$ are rotated $120^\circ$ clockwise and anti-clockwise, respectively)}
	\label{fig-VAcvged}
\end{figure}

A feasible solution is reported for NTUPF in Algorithm \ref{alg-1}. For comparison purpose, the problem is also solved when formulated in linear form, i.e. LTUPF. VMs and VAs of the simulation results are presented in Fig.\ref{fig-VMcvged} and Fig.\ref{fig-VAcvged}. 

Obviously, the errors of VMs and VAs are within acceptable limits demonstrating that the LTUPF is accurate enough to approximate the original non-convex problem. The calculated VAs also show that their values fall in $[-0.75^\circ,0.34^\circ],[119.29^\circ,120.39^\circ],[-120.71^\circ,119.61^\circ]$ for phase $a$, $b$ and $c$, respectively, which strongly support the assumption that VA of all nodes in each phase are with limited difference. The results also imply that the specified range for VA can be narrower to achieve higher accuracy of fitting the parameters in \eqref{fxy}.

\subsection{Case II}
Case II is reported infeasible by NTUPF, i.e. the algorithm is not converged, firstly when phase positions of PSD are set as initial values, implying the fundamental laws of providing electricity via the network cannot be satisfied. However, Algorithm \ref{alg-1} finally provides a feasible solution after switching the customers $33,41,57$ and $73$ from phase $b,b,c,c$ to phase $a,a,a,b$, respectively, demonstrating the efficacy of operational feasibility improvement brought by PSD.

\begin{figure}[H]
	\centering\includegraphics[scale=0.17]{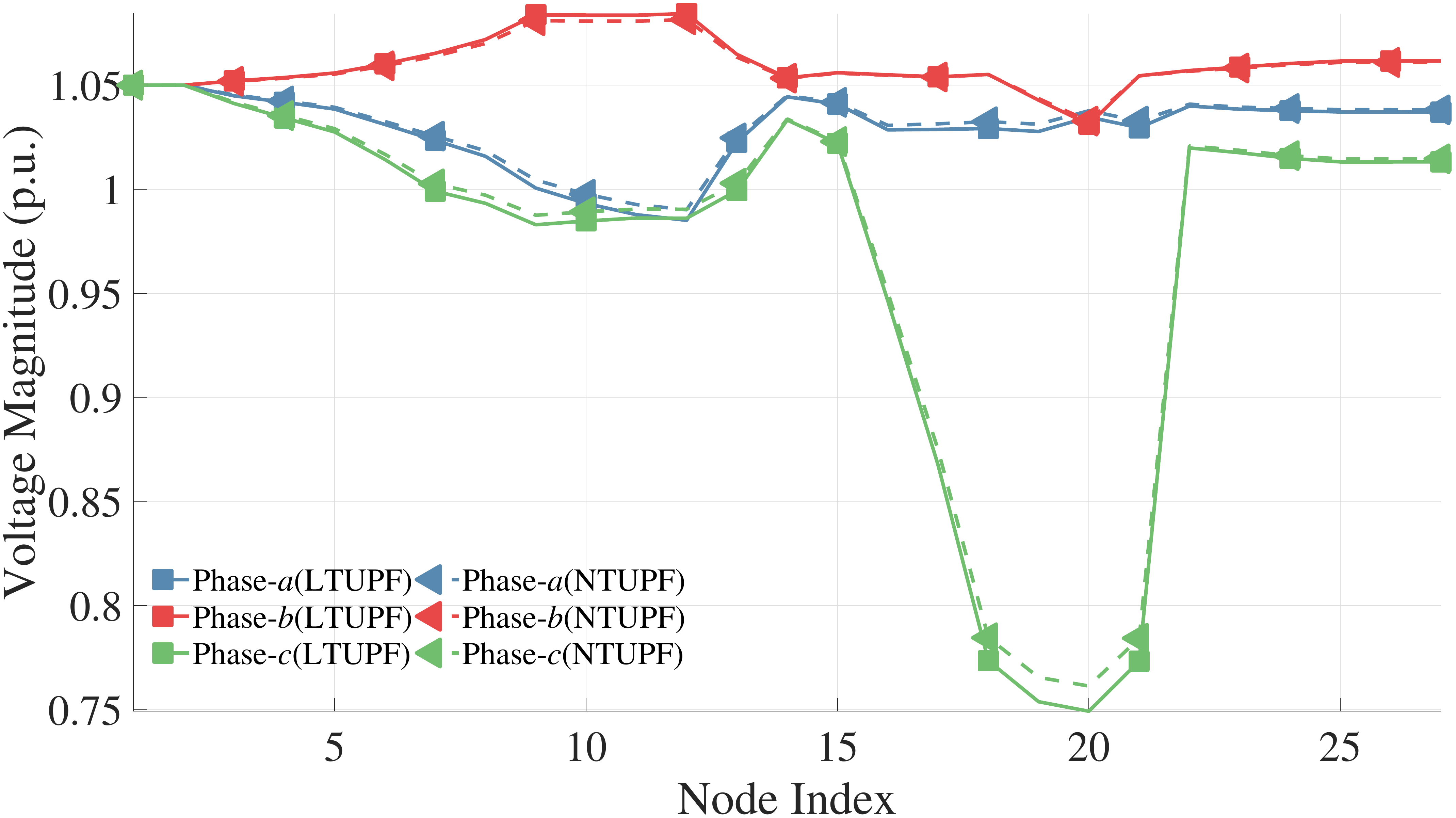}
	\caption{Calculated VMs for Case II}
	\label{fig-VMincvged}
\end{figure}

\begin{figure}[H]
	\centering\includegraphics[scale=0.17]{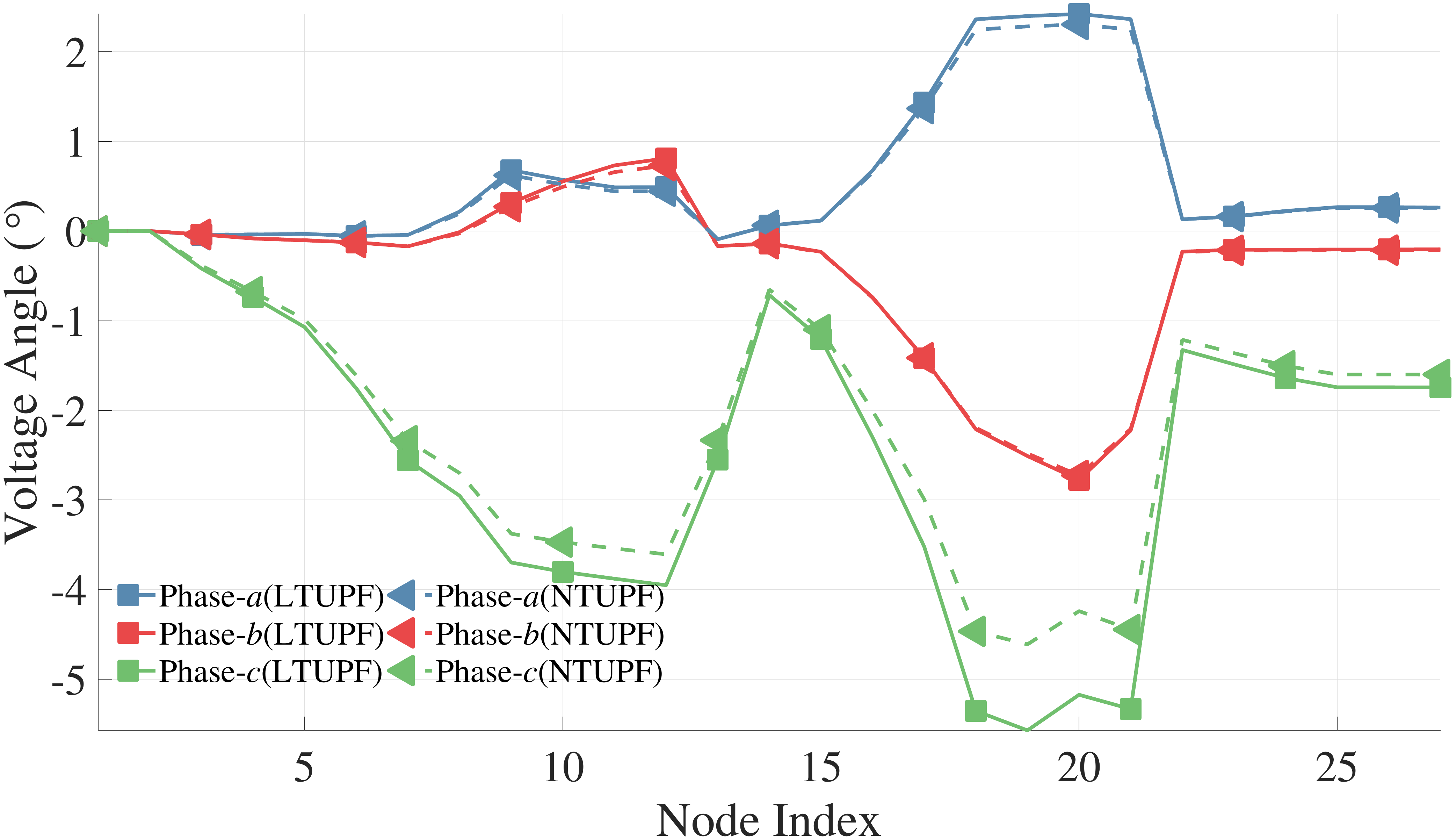}
	\caption{Calculated VAs for Case II (VAs for phase $b$ and $c$ are rotated $120^\circ$ clockwise and anti-clockwise, respectively)}
	\label{fig-VAincvged}
\end{figure}

Errors of VMs and VAs are presented in Fig.\ref{fig-VMincvged} and Fig.\ref{fig-VAincvged}, respectively, which again demonstrate that accuracy of approximating NTUPF by LTUPF is acceptable in engineering application.


\section{Conclusions}
Unbalance issue, which can be worsened by penetrating residential PV generation, may undermine the operational feasibility of LVDN, i.e. leading to infeasible TUPF problem. This paper presents an algorithm to switch the positions of PSDs to improve the operational feasibility and demonstrates its efficacy under some extreme load profiles. The proposed linear TUPF model for LVDN with PSD can be further developed to make it suitable in optimisation problems. Accordingly, other topics in this area, e.g. find the best locations to install PSDs, seeking optimal strategies to switch PSDs, all fall in our future research interests.

\section*{Acknowledgments}
The corresponding author of this paper would like to acknowledge that this work is funded by ARENA (Australian Renewable Energy Agency) Project: Demonstration of three-phase dynamic grid-side technologies for increasing distribution network DER hosting capacity. Disclaimer: The views expressed herein are not necessarily the views of the Australian Government, and the Australian Government does not accept responsibility for any information or advice contained herein.

\begin{spacing}{0.07}{\footnotesize
\bibliography{LPSDREF}

\begin{thebibliography}{10}
\providecommand{\url}[1]{#1}
\csname url@samestyle\endcsname
\providecommand{\newblock}{\relax}
\providecommand{\bibinfo}[2]{#2}
\providecommand{\BIBentrySTDinterwordspacing}{\spaceskip=0pt\relax}
\providecommand{\BIBentryALTinterwordstretchfactor}{4}
\providecommand{\BIBentryALTinterwordspacing}{\spaceskip=\fontdimen2\font plus
\BIBentryALTinterwordstretchfactor\fontdimen3\font minus
  \fontdimen4\font\relax}
\providecommand{\BIBforeignlanguage}[2]{{%
\expandafter\ifx\csname l@#1\endcsname\relax
\typeout{** WARNING: IEEEtran.bst: No hyphenation pattern has been}%
\typeout{** loaded for the language `#1'. Using the pattern for}%
\typeout{** the default language instead.}%
\else
\language=\csname l@#1\endcsname
\fi
#2}}
\providecommand{\BIBdecl}{\relax}
\BIBdecl

\bibitem{RN192}
\BIBentryALTinterwordspacing
``Solar report ({January} 2019),'' Australian Energy Council, Report, 2019.
  [Online]. Available:
  \url{https://www.energycouncil.com.au/media/15358/australian-energy-council-solar-report_-january-2019.pdf}
\BIBentrySTDinterwordspacing

\bibitem{RN61}
P.~K.~C. Wong, A.~Kalam, and R.~Barr, ``Modelling and analysis of practical
  options to improve the hosting capacity of low voltage networks for embedded
  photo-voltaic generation,'' \emph{IET Renew. Power Gener.}, vol.~11, no.~5,
  pp. 625--632, 2017.

\bibitem{RN52}
N.~K. Roy, H.~R. Pota, and M.~A. Mahmud, ``{DG} integration issues in
  unbalanced multi-phase distribution networks,'' in \emph{Proc. Australian
  Universities Power Engineering Conference ({AUPEC})}, 2016.

\bibitem{RN103}
J.~Zhu, M.-Y. Chow, and F.~Zhang, ``Phase balancing using mixed-integer
  programming,'' \emph{IEEE Trans. Power Syst.}, vol.~13, no.~4, pp.
  1487--1492, 1998.

\bibitem{RN104}
J.~Zhu, G.~Bilbro, and M.-Y. Chow, ``Phase balancing using simulated
  annealing,'' \emph{IEEE Trans. Power Syst.}, vol.~14, no.~4, pp. 1508--1513,
  1999.

\bibitem{RN44}
J.~Horta, D.~Kofman, D.~Menga, and M.~Caujolle, ``Augmenting der hosting
  capacity of distribution grids through local energy markets and dynamic phase
  switching,'' in \emph{Proc. the Ninth International Conference on Future
  Energy Systems - e-Energy '18}, Karlsruhe, Germany, 2018.

\bibitem{RN25}
T.~H. Chen, M.~S. Chen, K.~J. Hwang, P.~Kotas, and E.~A. Chebli, ``Distribution
  system power flow analysis-a rigid approach,'' \emph{IEEE Trans. Power Del.},
  vol.~6, no.~3, pp. 1146--1152, 1991.

\bibitem{RN26}
R.~D. Zimmerman and H.~D. Chiang, ``Fast decoupled power-flow for unbalanced
  radial-distribution systems,'' \emph{IEEE Trans. Power Syst.}, vol.~10,
  no.~4, pp. 2045--2052, 1995.

\bibitem{RN31}
H.~L. NGUYEN, ``Newton-raphson method in complex form,'' \emph{IEEE Trans.
  Power Syst.}, vol.~12, no.~3, pp. 1355--1359, 1997.

\bibitem{RN69}
P.~Garcia, J.~Pereira, S.~Carneiro, V.~d. Costa, and N.~Martins, ``Three-phase
  power flow calculations using the current injection method,'' \emph{IEEE
  Trans. Power Syst.}, vol.~15, no.~2, pp. 508 -- 514, 2000.

\bibitem{RN37}
T.~Jen-Hao, ``A direct approach for distribution system load flow solutions,''
  \emph{IEEE Trans. Power Del.}, vol.~18, no.~3, pp. 882--887, 2003.

\bibitem{RN70}
P.~Arboleya, C.~Gonzalez-Moran, and M.~Coto, ``Unbalanced power flow in
  distribution systems with embedded transformers using the complex theory in
  $\alpha\beta 0$ stationary reference frame,'' \emph{IEEE Trans. Power Syst.},
  vol.~29, no.~3, pp. 1012--1022, 2014.

\bibitem{RN55}
C.~Wang, A.~Bernstein, J.-Y. Le~Boudec, and M.~Paolone, ``Existence and
  uniqueness of load-flow solutions in three-phase distribution networks,''
  \emph{IEEE Trans. Power Syst.}, vol.~32, no.~4, pp. 3319--3320, 2018.

\bibitem{RN59}
A.~Bernstein, C.~Wang, E.~Dall'Anese, J.-Y. Le~Boudec, and C.~Zhao, ``Load flow
  in multiphase distribution networks: Existence, uniqueness, non-singularity
  and linear models,'' \emph{IEEE Trans. Power Syst.}, vol.~33, no.~6, pp.
  5832--5843, 2018.

\bibitem{RN23}
M.~Bazrafshan and N.~Gatsis, ``Convergence of the z-bus method for three-phase
  distribution load-flow with zip loads,'' \emph{IEEE Trans. Power Syst.},
  vol.~33, no.~1, pp. 153--165, 2018.

\bibitem{RN22}
------, ``Comprehensive modeling of three-phase distribution systems via the
  bus admittance matrix,'' \emph{IEEE Trans. Power Syst.}, vol.~33, no.~2, pp.
  2015--2029, 2018.

\bibitem{RN38}
H.~Ahmadi, J.~R. Marti, and A.~von Meier, ``A linear power flow formulation for
  three-phase distribution systems,'' \emph{IEEE Trans. Power Syst.}, vol.~31,
  no.~6, pp. 5012--5021, 2016.

\bibitem{RN39}
A.~Garces, ``A linear three-phase load flow for power distribution systems,''
  \emph{IEEE Trans. Power Syst.}, vol.~31, no.~1, pp. 827--828, 2016.

\bibitem{RN29}
A.~Bernstein and E.~Dall'Anese, ``Linear power flow models in multiphase
  distribution networks,'' in \emph{Proc. IEEE PES Innovative Smart Grid
  Technologies Conference Europe (ISGT-Europe)}, Turin, Italy, 2017.

\bibitem{RN30}
J.~A. Castrillon, J.~S. Giraldo, and C.~A. Castro, ``{MILP} for optimal
  reactive compensation and voltage control of distribution power systems,'' in
  \emph{Proc. IEEE Power \& Energy Society General Meeting}, Chicago, USA,
  2017.

\bibitem{BLiu-ITES}
B.~Liu, K.~Meng, Z.~Y. Dong, and W.~Wei, ``Optimal dispatch of coupled
  electricity and heat system with indepdent thermal energy storage system,''
  \emph{IEEE Trans. Power Syst.}, 2019.

\bibitem{RN50}
L.~Gan and S.~H. Low, ``Convex relaxations and linear approximation for opf in
  multiphase radial networks,'' in \emph{Proc. Power Systems Computation
  Conference}, Dublin, Ireland, 2018.

\bibitem{RN66}
B.~A. Robbins and A.~D. Dominguez-Garcia, ``Optimal reactive power dispatch for
  voltage regulation in unbalanced distribution systems,'' \emph{IEEE Trans.
  Power Syst.}, vol.~31, no.~4, pp. 2903--2913, 2016.

\bibitem{RN67}
J.~C. Lopez, J.~F. Franco, M.~J. Rider, and R.~Romero, ``Optimal
  restoration/maintenance switching sequence of unbalanced three-phase
  distribution systems,'' \emph{IEEE Trans. Smart Grid}, vol.~9, no.~6, pp.
  6058--6068, 2018.

\bibitem{BLiu-Bcoff}
B.~{Liu}, W.~{Wei}, and F.~{Liu}, ``Locating all real solutions of power flow
  equations: a convex optimisation-based method,'' \emph{IET Generation,
  Transmission Distribution}, vol.~12, no.~10, pp. 2273--2279, 2018.

\bibitem{gurobi}
\BIBentryALTinterwordspacing
L.~Gurobi~Optimization, ``Gurobi optimizer reference manual,'' 2018. [Online].
  Available: \url{http://www.gurobi.com}
\BIBentrySTDinterwordspacing

\bibitem{BPF-2}
B.~Liu, F.~Liu, B.~Zhai, and H.~Lan, ``Investigating continuous power flow
  solutions of {IEEE} 14-bus system,'' \emph{{IEEJ} Transactions on Electrical
  and Electronic Engineering}, vol.~14, no.~1, pp. 157--159, 2019.

\end{thebibliography}
}\end{spacing}
\end{multicols}

\end{document}